\documentclass[12pt]{amsart}

\usepackage{amssymb}

\usepackage{
setspace, subfig, caption}
\usepackage{url}
\usepackage{mathtools}

\usepackage{amssymb,amscd} 
\usepackage{enumerate} 
\usepackage[active]{srcltx} 
\usepackage{
leftidx}

\usepackage{multirow, booktabs}
\usepackage{algorithm}
\usepackage{algpseudocode}

\usepackage{graphicx}

\usepackage{tikz}
\usetikzlibrary{mindmap,shadows, calc}
\usepackage[hidelinks,pdfencoding=auto]{hyperref}

\newcommand{\mo}[1]{{
#1}}
\newcommand{\ds}[1]{{
#1}}

\usepackage[square,sort,comma,numbers]{natbib} 
\usepackage{algorithm}
\usepackage{algpseudocode}

\usepackage{hyperref}

\newcommand{\low}[1]{{\underline{#1}}}
\newcommand{\up}[1]{{\overline{#1}}}

\begin{document}
\numberwithin{equation}{section}

\theoremstyle{plain}
\newtheorem{proposition}{Proposition}[section]
\newtheorem{theorem}[proposition]{Theorem}
\newtheorem{lemma}[proposition]{Lemma}
\newtheorem{corollary}[proposition]{Corollary}

\theoremstyle{definition}
\newtheorem{definition}[proposition]{Definition}
\newtheorem{example}{Example}[section]
\newtheorem{remark}{Remark}[section]


\title[Maximal copulas]{A complete characterization of maximal copulas with given track section}
\author{Matja\v{z} Omladi\v{c}${}^{1}$}%
\thanks{${}^1$Supported by the ARIS (Slovenian Research and Innovation Agency) research core funding No. P1-0448.}
\address{Matja\v z Omladi\v c,
University of Ljubljana, Faculty of Mathematics and Physics, \& Institute of Mathematics, Physics and Mechanics, Jadranska 19, SI-1000 Ljubljana, Slovenia.} 
\email{matjaz@omladic.net}
\author{ Damjan \v{S}kulj${}^{2}$}
\thanks{${}^2$Supported by the ARIS (Slovenian Research and Innovation Agency) research core funding No. \ds{P5-0168}.}
\address{Damjan \v{S}kulj, University of Ljubljana, Faculty of Social Sciences\\
Kardeljeva pl. 5, SI-1000 Ljubljana, Slovenia}
\email{Damjan.Skulj@fdv.uni-lj.si}

\begin{abstract}
\ds{Bivariate copulas with prescribed diagonal section were first studied by Bertino. Their maximality was studied so far only from the point of view of upper bounds which brings quasi-copulas into the picture and limits the resulting set substantially. 
We propose to study maximality of these families in the order theoretic sense. A copula $C$ with given diagonal section $\delta$ is called \emph{undominated} if there is no copula $C' \neq C$ with the same diagonal section $\delta$ such that $C \leq C'$. The main contribution of this paper is a new method that provides copulas of the kind. 
This generates a much wider class that contains the known upper bounds as a very small subclass. There was a recent call for the study of asymmetry which is addressed by our class better than by the known ones.}
\mo{
Corresponding quasi-copulas can be obtained from our copulas via splicing techniques.
Most results are given on the level of tracks. } 
\end{abstract}
\maketitle

\section{Introduction}

It has been a problem to statisticians for many years how to construct a distribution with given marginals, until the Sklar’s theorem \cite{Skla} reduced it to the construction of a copula. We can think of copulas as distributions with uniform marginals and of quasi-copulas as pointwise suprema or infima of families of copulas, the exact definitions will be given in Section \ref{s-prlm}. In this paper we only consider bivariate copulas and quasi-copulas, so the word “bivariate” will be omitted.

Bivariate copulas with prescribed diagonal section were first studied by Bertino \cite{Bert}, and later by Fredricks and Nelsen \cite{FrNe}. The exact lower bound for such copulas is given in \cite{FrNe1} and it is the Bertino copula. The problem of upper bounds was settled in Nelsen et al.\ \cite{NeQuMoRoLaUbFl} and Kokol Bukov\v sek et al.\ \cite{KoBuMoSt}. It seems that the question was first set for tracks in Zou et al. \cite{ZoSuXi} and the result of \cite{KoBuMoSt} was extended to the track situation in Ouyang et al.\ \cite{OuSuZh}. (An interested reader may get acquainted with bivariate tracks in \cite{FeSaUbFl}.) We will unify the terminology by calling the two approaches the identity track case compared to the general track case. In all the papers cited so far the only view on maximality of copulas with given diagonal section goes through bounds; and the only other construction is what Nelsen et al.\ \cite{NeQuMoRoLaUbFl} are calling splicing -- one simply glues together the copulas with given diagonal obtained above and below the diagonal. In either way one obtains at most one quasi-copula of the desired type which is rarely a copula. 

We propose to study maximal copulas with given diagonal section in the order theoretic sense of the word. Concretely, a copula $C$ with given diagonal section $\delta$ is called \emph{undominated} if there is no copula $C' \neq C$ with the same diagonal section $\delta$ such that $C \leq C'$ with respect to the pointwise order. The main contribution of this paper is a new method that provides copulas of the kind. One of our main motivations to search for this method was the trend-setting paper of De Baets et al.\ \cite{DeBaDeMeJw} (cf.\ also Jwaid et al.\ \cite{JwDeMeHaIsDeBa}). They point out the 
(a)symmetry of a (quasi-)copula which is strongly related to that of exchangeability of random variables, a profoundly studied topic in probability theory and statistics. Symmetry of a function on the unit square is related to diagonal reflection with the main diagonal of the unit square acting as a mirror. In view of that our main result is Theorem \ref{thm-undominated} saying that for any distinct undominated copulas $C_1$ and $C_2$, i.e., $C_1\neq C_2$ with the same diagonal section, there exist points $u,v\in[0,1], u\neq v,$ such that
\[
    C_1(u,v)>C_2(u,v)\quad\mbox{and}\quad C_2(v,u)>C_1(v,u).
\] 
In other words, for two distinct undominated copulas with equal diagonal sections there exist mirror points with respect to the main diagonal such that each of them dominates the other in one of the two points.

The rest of the paper is structured as follows.
Section \ref{s-prlm} is devoted to preliminaries, Section \ref{s-mcgts} to the definition of our copulas. We start by a general track defined by the track function $\varphi$ and a $\varphi$-diagonal $\delta$. First assume we have a copula $C$ satisfying these data. Then we define four functions of one variable (either $x$ or $y$ from $[0,1]$) in geometric terms of this copula. These functions can all be uniquely determined from $\varphi$, $\delta$ and the first one of them, denoted by $\psi$. The canonical functions are always increasing and give rise to a function of two variables $C_\psi(x,y)$ which is always a copula and dominates $C$, i.e., $C\leq C_\psi$. 
Alternatively, we can start with $\varphi$, $\delta$ and $\psi$, use the previously established relations to define the rest of canonical functions and $C_\psi(x,y)$. If the canonical functions so defined are increasing, then $C_\psi$ is a copula and every undominated copula can be obtained in this way. Section \ref{sec-undom} is devoted to the study of the region, outside which copula $C_\psi$ always equals the Fr\' echet-Hoeffdung copula $M$. We also write these copulas in terms of variations, so that \emph{(a)} we can introduce a linear order on them, and \emph{(b)} we can find the known upper bounds of \cite{KoBuMoSt} among them. It turns out that only very rare of our undominated copulas have been found before. 
Section \ref{sec-final} presents our main results including Theorem \ref{thm-bounds-track} which exhibits equivalent definitions of copulas of the form $C_\psi$. Of course, once we know undominated copulas, we can get undominated quasi-copulas using the splicing techniques of \cite{JwDeMeHaIsDeBa}.

\section{Preliminaries}\label{s-prlm}

\mo{ Throughout this paper, we use the notation \( a \wedge b := \min\{a, b\} \).  
}

A function $C: [0,1]\times[0,1] \to [0,1]$ is called a \emph{quasi-copula} if conditions \emph{(a), (b),} and \emph{(c)} below are fulfilled. A \emph{copula} may be seen as a quasi-copula $C$ such that $C$-volume of every rectangle is positive.
\emph{\begin{enumerate}[(a)]
  \item $C$ is grounded and 1 is the {neutral element} of $C$;
  \item $C$ is increasing in each of its arguments on $[0,1]\times[0,1]$;
  \item $C$ satisfies the 1-Lipschitz condition on $[0,1]\times[0,1]$ in each of its arguments, i.e., it satisfies Condition (L):
\end{enumerate}}\noindent
\emph{\begin{enumerate}[(L)]
  \item \begin{enumerate}[(i)]
          \item $|C(x_2,b)-C(x_1,b)|\leq |x_2-x_1|$ for all $x_1,x_2,b \in [0, 1]$;
          \item $|C(a,y_2)-C(a,y_1)|\leq |y_2-y_1|$ for all $y_1,y_2,a \in [0, 1]$.
        \end{enumerate}
\end{enumerate}}
Here, condition \emph{(a)} means that $C(x, 0) = C(0,x)=0$ for all $x \in [0, 1]$ and $C(x, 1) = C(1,x)= x$ for all $x \in [0, 1]$. Condition \emph{(b)} means that\noindent \vskip-3mm

\emph{\noindent\begin{enumerate}[(I)]
  \item  $C(x_1,b)\leq C(x_2,b)$ for all $x_1,x_2,b \in [0, 1]$ such that $x_1\leq x_2$ and $C(a,y_1) \leq C(a,y_2)$ for all $y_1,y_2,a \in [0, 1]$ such that $y_1\leq y_2$.
\end{enumerate}}
Observe that we are using the term \emph{increasing} instead of more common ``non-decreasing''. The condition obtained from the first requirement of \emph{(I)} by replacing the two symbols $\leq$ by $<$ will be called \emph{strictly increasing (in $x$)} instead of more common ``increasing''. It remains to explain what we mean by ``$C$-volume of every rectangle is positive'': we mean that
\emph{\begin{enumerate}[(P)]
        \item $V_C([x_1, x_2] \times [y_1,y_2]) = C(x_2,y_2) - C(x_1, y_2) - C(x_2,y_1)+ C(x_1, y_1) \geq 0$  for all $x_1,x_2,y_1,y_2 \in [0, 1]$  such that $x_1\leq x_2$ and $y_1\leq y_2$.
      \end{enumerate}}
Here, a set of the form $[x_1, x_2] \times [y_1,y_2]\subseteq[0,1]\times[0,1]$, $x_1\leq x_2,y_1\leq y_2$, is called a \emph{rectangle} and condition \emph{(P)} may be seen as the definition of its $C$-volume denoted by $V_C$, together with the requirement that it is positive.

In this paper we are interested primarily in copulas which may be defined in a simpler although equivalent way. A \emph{copula} is a function $C: [0,1]\times[0,1] \to [0,1]$ satisfying conditions \emph{(a}) and \emph{(P)}. In the bivariate case we consider here, quasi-copulas have a similar definition. Namely, a \emph{quasi-copula} is a function $C: [0,1]\times[0,1] \to [0,1]$ satisfying conditions \emph{(a}) and
\emph{\begin{enumerate}[(P$\,'$)]
        \item $V_C(R) \geq 0$  for all rectangles $R\subseteq[0,1]\times[0,1]$ with nonempty intersection with the boundary of $[0,1]\times[0,1]$.
      \end{enumerate}}
\noindent(It turns out that in the multivariate case, this geometric definition of copula translates literally. However, the corresponding definition of quasi-copula becomes much more involved, cf.\ \cite[Section 2]{OmSt4}).\\[1mm]

Diagonal sections of bivariate copulas and partly quasi-copulas are the main object of our interest in this paper. A beautiful introduction into this subject may be found in \cite{DuKoMeSe}. The standard notation for a \emph{diagonal section} of a copula $C$ is $\delta$; we define $\delta(x)=C(x,x)$ for $x\in[0,1]$. It is well-known that
\begin{enumerate}[(a)]
  \item $\delta(x)\leq x$,
  \item $\delta(1)=1$,
  \item $\delta$ is 2-Lipschitz on $[0,1]$.
\end{enumerate}
Also, for any $\delta$ with these properties there is a copula having it as its diagonal section. The same is true for quasi-copulas.
\\[1mm]

    Now, we present some standard notions from real analysis as may be found, say, in \cite[Chapter 1]{Loja}. 
	For a function \(g: [0, 1] \to [0, 1]\), we define:
	\begin{enumerate}[(A)]
	  \item \emph{Total Variation on \([x, y]\)}: $TV_{x, y}(g) = \sup \sum_{i=1}^{n} |g(x_i) - g(x_{i-1})|$;
	  \item \emph{Positive Variation on \([x, y]\)}: 
\ds{\[ 
    V^+_{x, y}(g) = \sup \sum_{i=1}^{n} \max\left\{(g(x_i) - g(x_{i-1})) , 0\right\};
\]}
	  \item \emph{Negative Variation on \([x, y]\)}: 
\ds{ \[ 
    V^-_{x, y}(g) = \sup \sum_{i=1}^{n} \max\left\{(g(x_{i-1}) - g(x_i)) , 0\right\};
\]}
	\end{enumerate}
where the suprema are taken over all partitions \(x = x_0 < x_1 < \cdots < x_n = y\) of \([x, y]\) and all positive integers \(n\). If \(x = y\), set \(TV_{x, x}(g) = V^+_{x, x}(g) = V^-_{x, x}(g) = 0\).

We exhibit here, without proof, some standard properties for \(x \leq y\) needed in the sequel.
  \begin{enumerate}[(a)]
    \item $V^+_{x, y}(f) + V^-_{x, y}(f) = TV_{x, y}(f),$
    \item $V^+_{x, y}(f) - V^-_{x, y}(f) = f(y) - f(x),$
    \item $V^+_{x, y}(f) = V^+_{0, y}(f) - V^+_{0, x}(f),$
    \item $V^-_{x, y}(f) = V^-_{0, y}(f) - V^-_{0, x}(f),$
    \item $V^+_{x, y}(-f) = V^-_{x, y}(f).$
  \end{enumerate}
From {(a)} and {(b)}, we derive:
\begin{enumerate}[{(f)}]
	\item $V^+_{x, y}(f) = \frac{1}{2} \left( TV_{x, y}(f) + f(y) - f(x) \right),$
\end{enumerate}
\begin{enumerate}[{(g)}]
	\item $V^-_{x, y}(f) = \frac{1}{2} \left( TV_{x, y}(f) + f(x) - f(y) \right).$ 
\end{enumerate}

The following lemma is closely related to the Jordan canonical form (cf.~\cite[Section~1.3]{Loja}). We give an independent proof for the sake of completeness.

\begin{lemma}\label{lem-posvar}
  Let \( g: [0, 1] \to [0, 1] \) have bounded total variation. The function \( f(x) = V^+_{0, x}(g) \) is the minimal increasing function such that $f(0)=0$ and \( f - g \) is increasing. Moreover, every function \( h \) with the latter property satisfies \( h(y) - h(x) \geq V^+_{x, y}(g) \) for \( x \leq y \).
\end{lemma}

\begin{proof}
	First, note that \( f(x) = V^+_{0, x}(g) \) is increasing, since it is the supremum of sums of positive increments of \( g \). Moreover, \( f - g \) is increasing because
	\[
	f(y) - f(x) = V^+_{0, y}(g) - V^+_{0, x}(g) = V^+_{x, y}(g) \geq (g(y) - g(x)) \vee 0 \geq g(y) - g(x).
	\]
	
	To prove minimality, let \( h \) be any increasing function such that \( h - g \) is increasing, i.e., \( h(v) - h(u) \geq g(v) - g(u) \) for all \( u \leq v \). Given \( \varepsilon > 0 \), there exists a partition \( 0 = x_0 < x_1 < \cdots < x_n = x \) such that
	\[
	\sum_{i=1}^{n} (g(x_i) - g(x_{i-1})) \vee 0 > V^+_{0, x}(g) - \varepsilon.
	\]
	Using monotonicity of \( h \), we get
	\[
	h(x) - h(0) = \sum (h(x_i) - h(x_{i-1})) \geq \sum (g(x_i) - g(x_{i-1})) \vee 0 > V^+_{0, x}(g) - \varepsilon.
	\]
	If \( g(0) = 0 \)
, then \( h(0) \geq g(0) = 0 \), so \( h(x) \geq V^+_{0, x}(g) \), proving minimality.
	
	Finally, since \( h(y) - h(x) \geq g(y) - g(x) \), and using the identity
	\[
	V^+_{x, y}(g) = V^+_{0, y}(g) - V^+_{0, x}(g),
	\]
	we conclude \( h(y) - h(x) \geq V^+_{x, y}(g) \).
\end{proof}

We conclude this section with a simple observation that a bivariate function constructed as the minimum of two increasing univariate functions produces positive volumes of rectangles. This well-known fact makes a copula out of the Fr\' echet-Hoeffding  upper bound of all bivariate copulas \( M(x, y) = x \wedge y \). Since this simple truth plays a key role in our investigations, we present its proof for the sake of completeness.
Observe that the value of the left-hand side of the inequality below is exactly the volume of the rectangle $[x_1,x_2]\times[y_1,y_2]$. 

\begin{proposition}\label{prop-maximum-positive-volumes}
	Let \ds{\( f , g
\) be increasing real functions defined each on an interval}, and define \( h(x, y) = f(x) \wedge g(y) \). Then we have for any \( x_1 < x_2 \) and \( y_1 < y_2 \) \ds{from respective domains of the functions}
	\begin{equation*}
		h(x_1, y_1) + h(x_2, y_2) - h(x_1, y_2) - h(x_2, y_1) \geq 0.
	\end{equation*}
\end{proposition}

\begin{proof}
	Denote the value of the left-hand side of this inequality by $V$ and let us prove it case by case.

	\noindent	\textbf{Case 1.} \( f(x_1) \leq g(y_1) \leq g(y_2) \) and \( f(x_2) \leq g(y_2) \)]  
		Here,  $h(x_1, y_1) = f(x_1)= h(x_1, y_2)$ so that $V= f(x_2) - h(x_2, y_1) \geq 0$.

	\noindent	\textbf{Case 2.} \( f(x_1) \leq g(y_1) \leq g(y_2) \leq f(x_2) \)]  
		This time $V= g(y_2) - g(y_1) \geq 0$.

	The remaining cases follow by symmetry, interchanging \( f \) and \( g \), thus ensuring non-negativity in all scenarios.
\end{proof}

\section{Undominated copulas for general tracks }\label{s-mcgts}

In this section we develop a new method that finds all (!) copulas that are maximal with a given diagonal on a track. 

Let us specify. First, we replace the term ``maximal'', which is used in many meanings even in copula theory, by a more suggestive and unambiguous term throughout the paper. We call a copula $C$ \emph{undominated} (in a certain family) if there is no copula $C' \neq C$ (in the same family) such that $C \leq C'$. 

Next, we use a strictly increasing, continuous function \( \varphi\colon [0, 1] \to [0, 1] \) with \( \varphi(0) = 0 \) and \( \varphi(1) = 1 \) to define a \emph{track}:
	\[
	B_\varphi = \{ (x, \varphi(x)) \mid  x \in [0, 1] \}.
	\]
Strict monotonicity ensures the existence of a unique inverse \( \varphi^{-1} \).  Note that track here is seen as the graph of function $\varphi$ rather than the function itself. \mo{We call $\varphi$ the\emph{ track function}.}
We follow the routine approach (cf. \cite{FeSaUbFl}) in this and the next definition. 

We say that the function $\delta : [0, 1] \rightarrow [0, 1]$ is a \emph{$\varphi$-diagonal} if it satisfies the following conditions:
\emph{\begin{enumerate}[(a)]  
  \item $\delta(1) = 1$,
  \item $\delta(t) \le \min \{t, \varphi(t)\}$ for all $t\in  [0, 1]$,
  \item $\delta$ is increasing, and
  \item $|\delta(t') - \delta (t)| \le |t - t'| + |\varphi(t') - \varphi(t)|$  for all $t,t'\in [0, 1]$.
  \end{enumerate}
}
	For a copula \( C \) and a track \( B_\varphi \), the function
	\[
	\delta(x) = C(x, \varphi(x))
	\]
	is called the \emph{track section} corresponding to \( C \) and \( B_\varphi \) (or simply corresponding to a track function $\varphi$). It is known \cite{FeSaUbFl} that a function $\delta : [0, 1] \rightarrow [0, 1]$ is a {track section} corresponding to a track function $\varphi$ if and only if it is a $\varphi$-diagonal. 

So, given a track function $\varphi$ and a  $\varphi$-diagonal $\delta$ we want to know all (!) undominated copulas with $\delta$ as a track section.

Let $\varphi$ be any track function according to the above definition and let $\delta$ be a $\varphi$-diagonal. Introduce
the set \( \mathcal{C}(B_\varphi, \delta) \) of all copulas \( C \) whose track section corresponding to $\varphi$ equals $\delta$. 
We define the following subsets of $[0, 1]^2$:
\begin{align*}
	L_\varphi & = \{ (u, v) \mid  v\le \varphi (u) \} &
	U_\varphi & = \{ (u, v) \mid  v>\varphi(u) \} \\
	M(x) & = \{ (u, v) \mid  u \le x  \} &
	N(y) & = \{ (u, v) \mid  v \le y  \} \\
	\Psi_\varphi(x) & = L_\varphi \cap M(x) &
	\Gamma_\varphi(y) & = L_\varphi \cap N (y) \\
	\Delta_\varphi(y) & = U_\varphi \cap N (y) &
	\Xi_\varphi(x) & = U_\varphi \cap M(x) 
\end{align*}
Every copula $C$ induces a probability measure $P_C$ on the family of Borel subsets of $[0, 1]^2$.

\begin{proposition}\label{prop-cont-initial} Fix a copula $C\in \mathcal C(B_\varphi, \delta)$ and introduce functions
\begin{align*}
	\psi(x) &= P_C(\Psi_\varphi(x)), &
	\chi(y) &= P_C(\Gamma_\varphi(y)), \\
	\eta(y) &= P_C(\Delta_\varphi(y)),  &
	\xi(x) & = P_C(\Xi_\varphi(x)).  
\end{align*}
Then, these functions are 
\begin{enumerate}[(a)]
  \item continuous, increasing on $[0,1]$, and
  \item satisfy the boundary conditions:
	\[
	\psi(0) = \chi(0) = \eta(0) = \xi(0) = 0, \quad \psi(1) = \chi(1) = 1-\eta(1) = 1-\xi(1).
	\]
\end{enumerate}
\end{proposition}

\begin{proof}
	Monotonicity follows directly from the definitions. Continuity is a consequence of the Lipschitz property of copulas:
	\[
	|C(x+h, y+k) - C(x, y)| \leq |h| + |k|.
	\]
	In particular, using $C(x, \varphi(x)) = \psi(x) + \eta(\varphi(x))$, we obtain, for any $h > 0$,
	\[
\begin{split}
   \psi(x+h) - \psi(x) & = C(x+h, \varphi(x+h)) - \eta(\varphi(x+h)) - C(x, \varphi(x)) + \eta(\varphi(x)) \\
   & \leq C(x+h, \varphi(x+h)) - C(x, \varphi(x)) \leq h + (\varphi(x+h) - \varphi(x)).
\end{split}
	\]
	Since $\varphi$ is continuous, the right-hand side tends to $0$ as $h \to 0$. Similar arguments apply to $\chi, \eta$, and $\xi$.
	
	The boundary conditions follow from the definitions and the copula properties. Specifically, $\psi(0) = P_C(\Psi_\varphi(0)) = 0$, and similarly for $\chi(0), \eta(0)$, and $\xi(0)$. At the upper boundary, we have
	\[
	\psi(1) = P_C(\Psi_\varphi(1)) = P_C(L_\varphi) = \chi(1),
	\]
	while 
	\[
	\eta(1) = P_C(U_\varphi) = 1 - P_C(L_\varphi) = 1 - \psi(1),
	\]
	and similarly, $\xi(1) = 1 - \psi(1)$. This completes the proof.
\end{proof}

The following proposition tells us that once we fixed a track section $\delta$ corresponding to a track function $\varphi$ the quadruplet of functions $(\psi, \chi, \eta, \xi)$ introduced in Proposition \ref{prop-cont-initial} is uniquely determined by its first member $\psi$. This will simplify our notation in the sequel.

\begin{proposition}\label{prop-relations}
	Let $C$ be a copula in $\mathcal{C}(B_\varphi, \delta)$. Then for the quadruplet $(\psi, \chi, \eta, \xi)$ of Proposition \ref{prop-cont-initial} the following relations hold:  
\begin{enumerate}[(a)]
  \item $\eta(y) = \delta(\varphi^{-1}(y)) - \psi(\varphi^{-1}(y)), $
  \item $\chi(y) = y - \eta(y) = y - \delta(\varphi^{-1}(y)) + \psi(\varphi^{-1}(y)), $
  \item $\xi(x) = x - \psi(x). $
\end{enumerate}
\end{proposition}

Any tuple of functions $(\psi, \chi, \eta, \xi)$ satisfying the conclusions of Proposition \ref{prop-cont-initial} and the relations of Proposition \ref{prop-relations} will be called \emph{canonical quadruplet}. As pointed out above, it is uniquely determined by the function $\psi$ (and, of course, by $\varphi$ and $\delta$).

\begin{proof}
	By definition of the probability measure induced by $C$, we have
	\[
	P_C(M(x) \cap N(\varphi(x))) = C(x, \varphi(x)) = \delta(x).
	\]
	Since 
	\[
	\Psi_\varphi(x) \cup \Delta_\varphi(\varphi(x)) = M(x) \cap N(\varphi(x))
	\]
	with $\Psi_\varphi(x) \cap \Delta_\varphi(\varphi(x)) = \emptyset$, it follows that
	\[
	\psi(x) + \eta(\varphi(x)) = \delta(x).
	\]
	Substituting $y = \varphi(x)$ and $x = \varphi^{-1}(y)$ the above relation gives
	\[
	\eta(y) = \delta(\varphi^{-1}(y)) - \psi(\varphi^{-1}(y)),
	\]
	which proves \emph{(a)}.
	
	Next, since  
	\[
	\Delta_\varphi(y) \cup \Gamma_\varphi(y) = N(y), \quad \Delta_\varphi(y) \cap \Gamma_\varphi(y) = \emptyset,
	\]
	we deduce that  
	\[
	\eta(y) + \chi(y) = P_C(N(y)) = y.
	\]
	Rearranging this equation and using \emph{(a)}, we obtain
	\[
	\chi(y) = y - \eta(y) = y - \delta(\varphi^{-1}(y)) + \psi(\varphi^{-1}(y)),
	\]
	proving \emph{(b)}.
	
	Finally, since  
	\[
	\Xi_\varphi(x) \cup \Psi_\varphi(x) = M(x), \quad \Xi_\varphi(x) \cap \Psi_\varphi(x) = \emptyset,
	\]
	it follows that  
	\[
	\xi(x) + \psi(x) = P_C(M(x)) = x.
	\]
	Rearranging gives  
	\[
	\xi(x) = x - \psi(x),
	\]
	which proves \emph{(c)} as well.
	\qedhere
\end{proof}

We now write any copula from $\mathcal{C}(B_\varphi, \delta)$ as a sum of two functions which will play a key role in our further investigations.

\begin{proposition}\label{prop-ineq-S-T}
	Let $C\in \mathcal C(B_\varphi, \delta)$ be a copula, and define 
\[
	S(x, y) = P_C(M(x) \cap N(y) \cap L_\varphi)\ \mbox{and}\ T(x, y) = P_C(M(x) \cap N(y) \cap U_\varphi); 
\]
then:
\begin{enumerate}[(a)]
  \item $S(x, y) + T(x, y) = C(x, y)$;
  \item $S(x, y) \leq \psi(x) \wedge \chi(y)$;
  \item $T(x, y) \leq \eta(y) \wedge \xi(x)$.
\end{enumerate}
	where $(\psi, \chi, \eta, \xi)$ is the canonical quadruplet of Proposition~\ref{prop-cont-initial}.
\end{proposition}

\begin{proof}
\mo{Since \( M(x) \cap N(y) = (M(x) \cap N(y) \cap L_\varphi) \cup (M(x) \cap N(y) \cap U_\varphi) \), we represented set $M(x) \cap N(y)$ as a union of two sets with empty intersection, so that \emph{(a)} follows. }	
From the inclusion $M(x) \cap N(y) \cap L_\varphi = \Psi_\varphi(x) \cap N(y) \subseteq \Psi_\varphi(x)$, we obtain $S(x, y) \leq \psi(x)$. Similarly, since $M(x) \cap N(y) \cap L_\varphi = \Gamma_\varphi(y) \cap M(x) \subseteq \Gamma_\varphi(y)$, it follows that $S(x, y) \leq \chi(y)$. Taking the minimum of these bounds proves \emph{(b)}.
	For $T(x, y)$, the inclusion $M(x) \cap N(y) \cap U_\varphi = \Delta_\varphi(y) \cap M(x) \subseteq \Delta_\varphi(y)$ implies $T(x, y) \leq \eta(y)$. Likewise, since $M(x) \cap N(y) \cap U_\varphi = \Xi_\varphi(x) \cap N(y) \subseteq \Xi_\varphi(x)$, we get $T(x, y) \leq \xi(x)$. Taking the minimum of these bounds gives \emph{(c)}.
\end{proof}

 Every copula $C\in \mathcal C(B_\varphi, \delta)$ induces a function $\psi$ and the related functions, such that inequalities from Proposition~\ref{prop-ineq-S-T} are satisfied. Next we show that these inequalities may become exact by constructing functions $S_\psi$ and $T_\psi$, which achieve the equalities in these relations. 

However, we will observe in the following theorem that this result is independent of the choice of the starting copula $C \in \mathcal{C}(B_\varphi, \delta)$. It will only depend on the choice of the function $\psi$ satisfying certain conditions. Choose a track function $\varphi$, a function $\delta$ satisfying all conditions of the definition for a $\varphi$-diagonal except possibly for condition \emph{(d)}\footnote{Condition \emph{(d)} will follow from the other assumptions.}, and an appropriate function $\psi$. Let
\emph{\begin{enumerate}[(a)]
  \item $S_\psi (x, y) = \psi(x) \wedge (y - \delta(\varphi^{-1}(y)) + \psi(\varphi^{-1}(y))),$
  \item $T_\psi (x, y) = (x - \psi(x)) \wedge (\delta(\varphi^{-1}(y)) - \psi(\varphi^{-1}(y))),$
  \item $C_\psi (x,y) = S_\psi (x,y) + T_\psi (x,y).$
\end{enumerate}	}
In order to simplify these definitions we use canonical functions defined using the relations of Proposition \ref{prop-relations}\emph{(a),(b),(c)}
to get
\begin{enumerate}[(a')]
  \item \ds{$S_\psi(x,y)=\psi(x)\wedge\chi(y),$}
  \item \ds{$T_\psi(x,y)=\xi(x)\wedge\eta(y),$}
  \item $C_\psi(x,y)=S_\psi(x,y)+T_\psi(x,y).$
\end{enumerate}

\begin{theorem}\label{prop-copula-psi} Choose a track function $\varphi$ and a function $\delta$ as above.
	Let $\psi:[0,1]\to[0,1]$ be increasing with $\psi(0)=0$ and define $S_\psi, T_\psi, C_\psi: [0,1]^2 \to [0,1]$ as above. Assume also that $\chi,\eta,\xi$ defined as above are non-negative and increasing on $[0,1]$.
	
	Then $C_\psi$ is a copula that matches the prescribed track section:
	\[
	C_\psi(x,\varphi(x))=\delta(x)\qquad \text{for all }x\in[0,1].
	\]
\end{theorem}


\begin{proof}
	By assumption all canonical functions are 
non-negative and increasing. Hence, by Proposition~\ref{prop-maximum-positive-volumes}, both $S_\psi$ and $T_\psi$ induce positive volumes on every rectangle $R=[x_1, x_2]\times [y_1, y_2]$. Since $C_\psi$ is their sum, it induces positive volumes as well. 
	
	To conclude $C_\psi$ is a copula, we need to verify its marginals. It is not hard to see that our canonical functions as defined above automatically satisfy conditions of Proposition \ref{prop-cont-initial}\emph{(b)}.
We next observe that $C_\psi(0, x) = C_\psi(x, 0) = 0$. Now take $x = 1$. Since $\chi(y)$ is increasing and $\chi(1)= \psi(1)$, it follows that $S_\psi(1, y) = y-\delta(\varphi^{-1}(y)) + \psi(\varphi^{-1}(y))$ and similarly, $T_\psi(1, y) = \delta(\varphi^{-1}(y)) - \psi(\varphi^{-1}(y))$. It follows that $C_\psi(1, y) = y$. Similarly, we get that $S_\psi(x, 1) = \psi(x)$ and $T_\psi(x, 1) = x-\psi(x)$, implying that $C_\psi(x, 1) = x$. 
	
	Finally, we calculate the value of $C_\psi(x, \varphi(x))$. Note that $\delta(x)  \le \varphi(x)\wedge x \le \varphi(x)$. 
This implies that $\psi(x) \le  \varphi(x)-\delta(x)+\psi(x)$ and therefore $S_\psi(x, \varphi(x)) = \psi(x)$. Similarly, because of $\delta(x)\le x$ we have that $T_\psi(x, \varphi(x)) = \delta(x)-\psi(x)$. Summing up we get $C_\psi(x, \varphi(x)) = \psi(x)+ (\delta(x)-\psi(x)) = \delta(x)$. 
\end{proof}


\begin{corollary}\label{cor-c-psi-undominated}
	Let $C$ be a copula satisfying $C(x, \varphi(x)) = \delta(x)$ for all $x \in [0,1]$. If $C$ is undominated, then there exists a function $\psi$, satisfying conditions of Theorem \ref{prop-copula-psi}, such that $C = C_\psi$. 
\end{corollary}
\begin{proof}
	Define $\psi$ as in Proposition \ref{prop-cont-initial}. By Proposition~\ref{prop-ineq-S-T}, we have $C \leq C_\psi$. Since $C$ is undominated, it follows that $C = C_\psi$. 
\end{proof}

Next we provide an alternative formula for $C_\psi$ based on an observation that it can differ from the Fr\'echet-Hoeffding upper bound only on a region bounded by two functions, $$g(x)  = \max\{ y\mid  \psi(\varphi^{-1}(y)) + y-\delta(\varphi^{-1}(y)) \le \psi(x) \}$$
and $$h(x)  = \max\{ y\mid  \xi(\varphi^{-1}(y)) + \delta(\varphi^{-1}(y)) - \varphi^{-1}(y) \le \xi(x)\}.$$ 
This region contains the track section. \mo{The shape of the region will be studied in Section \ref{sec-undom} and some figures will be given there.} 
Let us recall that $\xi(x) = x-\psi(x)$.

\begin{proposition}\label{prop-g-h}
	The following holds for functions $g$ and $h$:
	\begin{enumerate}[(a)]
		\item The functions $g$ and $h$ are increasing.
		\item We have $g(x) \leq \varphi(x) \leq h(x)$ for all $x\in[0,1]$.
		\item For all $x\in[0,1]$,
		\[
		\varphi(x)-\delta(x)>0 \quad\Longleftrightarrow\quad g(x)<\varphi(x).
		\]
		\item For all $x\in[0,1]$,
		\[
		x-\delta(x)>0 \quad\Longleftrightarrow\quad \varphi(x)<h(x).
		\]
	\end{enumerate}
\end{proposition}

\begin{proof}
	(a) Let $x'<x$. Since $\psi$ is increasing, $\psi(x')\geq\psi(x)$, so
	\begin{multline*}
	\{y\mid \psi(\varphi^{-1}(y)) + y - \delta(\varphi^{-1}(y)) \leq \psi(x)\} \\
	\subseteq
	\{y\mid \psi(\varphi^{-1}(y)) + y - \delta(\varphi^{-1}(y)) \leq \psi(x')\},
	\end{multline*}
	implying $g(x)\leq g(x')$. Similarly, since $\xi(x)=x-\psi(x)$ is increasing, $\xi(x') \ge \xi(x)$, and 
	\begin{multline*}
	\{y\mid \xi(\varphi^{-1}(y)) + \delta(\varphi^{-1}(y)) - \varphi^{-1}(y)\leq\xi(x)\} \\
	\subseteq 
	\{y\mid \xi(\varphi^{-1}(y)) + \delta(\varphi^{-1}(y)) - \varphi^{-1}(y)\leq\xi(x')\},
	\end{multline*}
	so $h(x)\leq h(x')$. Thus both $g$ and $h$ are increasing.
	
	(b) Taking $y=\varphi(x)$, we have $\varphi^{-1}(y)=x$ and thus
	\[
	\psi(\varphi^{-1}(\varphi(x)))+\varphi(x)-\delta(\varphi^{-1}(\varphi(x)))
	=\psi(x)+\varphi(x)-\delta(x)\geq\psi(x).
	\]
	This implies that $\varphi(x)$ is an upper bound for $g(x)$, which by monotonicity of $\psi+\varphi- \delta$ implies that $g(x)\leq\varphi(x)$.
	Similarly,
	\[
	\xi(\varphi^{-1}(\varphi(x)))+\delta(\varphi^{-1}(\varphi(x)))-\varphi^{-1}(\varphi(x))
	=\xi(x)+\delta(x)-x\leq\xi(x)
	\]
	since $x-\delta(x)\geq0$. So $\varphi(x)$ is also an upper bound for $h(x)$, which by monotonicity of $\xi(x)+\delta(x)-x = \delta(x)-\psi(x)$ implies $\varphi(x)\leq h(x)$.
	
	(c) If $\varphi(x)-\delta(x) > 0$ then taking $y = \varphi(x)$ yields the left hand side in the definition of function $g$ is strictly larger than the right hand side, whence by continuity of all involved functions, implies that $g(x)<\varphi(x)$; while the case $\varphi(x)-\delta(x) = 0$, yields $g(x) = \varphi(x)$. Hence, the equivalence follows.

	(d) follows by similar reasoning
.
\end{proof}


\begin{corollary}\label{cor-c-psi-expression-track}
	The copula $C_\psi$ can be given in terms of the following formula
\begin{align*}
	C_\psi(x, y) & = \begin{cases}
						x, & h(x) \le y \\
						\psi(x)-\psi(\varphi^{-1}(y)) + \delta(\varphi^{-1}(y)), & g(x) \le y \le h(x) \\
						y, & y \le g(x)
					\end{cases} \\
					& = \min\{ x, y, \psi(x)-\psi(\varphi^{-1}(y)) + \delta(\varphi^{-1}(y))\}. \label{eq-c-psi-track-as-min}
\end{align*}
\end{corollary}
\begin{proof}
	The formula follows directly from the definition of functions $g$ and $h$, and  \mo{Theorem \ref{prop-copula-psi} together with the definitions \emph{(a),(b),(c)} given in its preamble.}
\end{proof}

Denote
\[
    \mathcal{R}_\psi=\{(x,y)\mid g(x)\leq y\leq h(x)\}.
\]
Clearly, outside the region $\mathcal{R}_\psi$ copula $C_\psi$ equals the min-copula $M$, i.e., the Fr\'echet-Hoeffding upper bound, while inside the region it depends on $\psi$. 

In the next section we restrict the analysis to the case where the track is diagonal, $\varphi(x) = x$. There we provide some simplifications and some stronger results. 

\section{Undominated copulas }\label{sec-undom} 

In this section, we fix function $\varphi$ in the track \(B_\varphi = \{ (x, \varphi(x)) \mid x \in [0, 1] \}\) as the identity function, i.e., $\varphi(x)=x$ for all $x \in [0, 1]$. So, the track section becomes \emph{diagonal section} \(\delta(x) = C(x, x)\) of a copula \(C\) on \([0, 1]^2\). Function \(\delta(x)\) prescribes the copula’s values along the \emph{diagonal} \(\{ (x, x) \mid x \in [0, 1] \}\). To emphasize the difference between the general track version we will call this version the \emph{identity track version}.

Recall canonical functions $\psi$ and \(\xi(x) = x - \psi(x)\) from the previous section and introduce also \(\zeta(x) = x - \delta(x)\). Since $\delta(x)\leq x$ and by Proposition \ref{prop-cont-initial}\emph{(a)} also $\psi(x)\leq \delta(x)$, we conclude \(\psi(x) \leq x\). It follows that
both \(\zeta(x)\) and \(\xi(x)\) are non-negative functions and $\xi$ is also an increasing function on \([0, 1]\).

Before we rewrite in a simpler form (using functions $\xi$ and $\zeta$) our undominated copula $C_\psi$ together with the two functions $g$ and $h$ that determine the region $\mathcal{R}_\psi$ outside which it equals copula min, we make an observation about this region for the diagonal case.

\begin{proposition}\label{prop-delta-0-min}
	Let \(C\) be a copula with diagonal section \(\delta(x) = C(x, x)\). If \(\zeta(x) \wedge \zeta(y) = 0\) for some $x,y\in[0,1]$,  
then \(C(x, y) = x \wedge y\).
\end{proposition}
\begin{proof}
	Assume without loss of generality that \(x \leq y\). If \(\zeta(x) \wedge \zeta(y) = 0\), then either \(\zeta(x) = 0\) or \(\zeta(y) = 0\).
	
	\textbf{Case 1: \(\zeta(x) = 0\)}: Here, \(\delta(x) = x\), so \(C(x, x) = x\). By monotonicity, \(C(x, y) \geq C(x, x) = x\), and since \(C(x, y) \leq C(x, 1) = x\), we have \(x \geq C(x, y) \geq x\). Thus, \(C(x, y) = x = x \wedge y\).
	
	\textbf{Case 2: \(\zeta(y) = 0\)}: Here, \(\delta(y) = y\), so \(C(y, y) = y\). Using the Lipschitz property, \(|C(y, y) - C(x, y)| \leq |y - x|\), so \(y - C(x, y) \leq y - x\), implying \(C(x, y) \geq x\). Also, \(C(x, y) \leq C(x, 1) = x\). Hence, \(x \geq C(x, y) \geq x\), so \(C(x, y) = x = x \wedge y\).
\end{proof}

\begin{corollary}\label{cor-delta-t-min}
	If \(x \leq y\) and \(t \in [x, y]\) are such that \(\zeta(t) = 0\). Then \(C(x, y) = C(y, x) = x \wedge y = x\).
\end{corollary}
\begin{proof}
	Since \(\zeta(t) = 0\), Proposition~\ref{prop-delta-0-min} gives \(C(x, t) = x \wedge t = x\) (as \(x \leq t\)). By monotonicity, \(C(x, t) \leq C(x, y) \leq C(x, 1) = x\), so \(x \leq C(x, y) \leq x\), implying \(C(x, y) = x\).
	
	Similarly, \(C(t, x) = t \wedge x = x\) (as \(t \geq x\)). Since \(t \leq y\), monotonicity gives \(C(t, x) \leq C(y, x) \leq C(1, x) = x\), so \(x \leq C(y, x) \leq x\), implying \(C(y, x) = x\). Thus, \(C(x, y) = C(y, x) = x = x \wedge y\).
\end{proof}

\begin{center}
\begin{figure}[h!]
      \caption{Images corresponding to $\delta(x) = x- \sin^2(2\pi x)/(2\pi)$. We exhibit the image of copula (up) and region (down) for the lower case $\psi_L$ (left) and the upper one $\psi_U$ (right).}\label{fig-bimod}
\includegraphics[width=5cm]{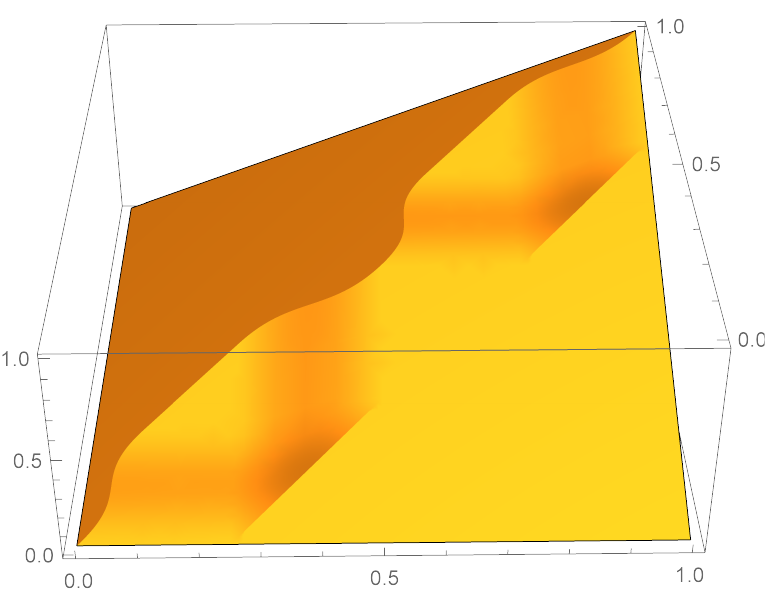}\ \ \includegraphics[width=5cm]{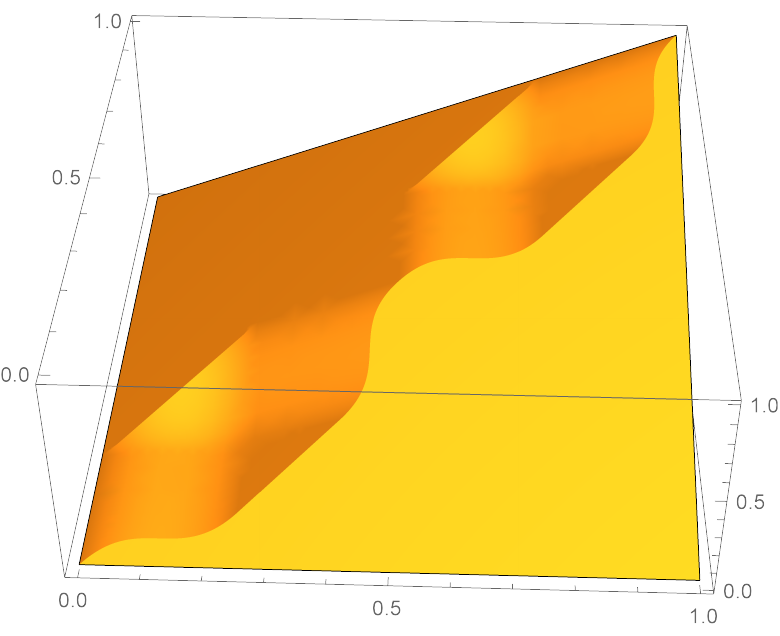}\ \ \includegraphics[width=4cm]{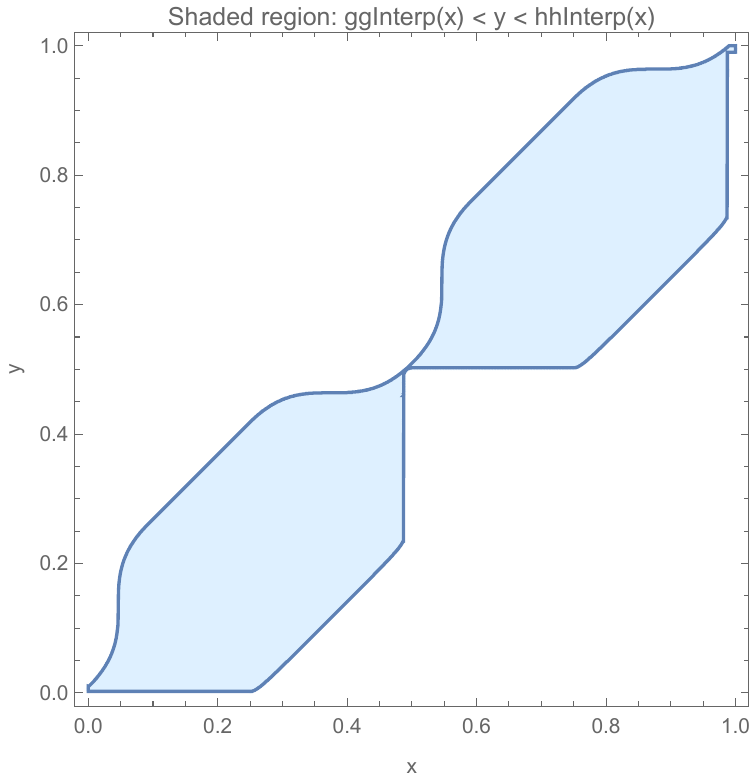}\ \ \ \ \ \ \includegraphics[width=4cm]{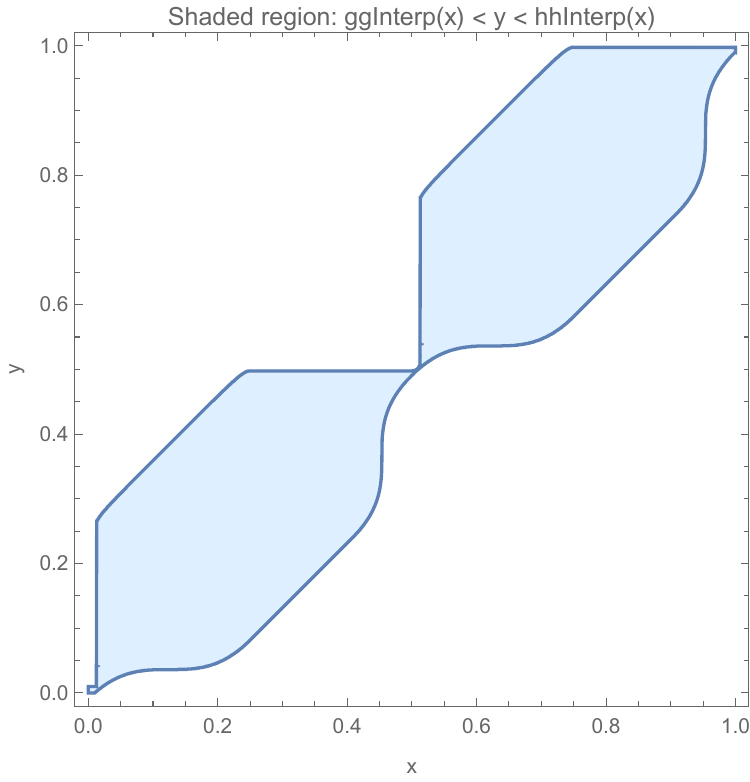}
\label{fig:figure1}
\end{figure}
\end{center}

\begin{center}
\begin{figure}[h!]
      \caption{Images corresponding to $\delta(x) = x-\sin(\pi x)/\pi$. We exhibit the image of copula (up) and region (down) for the lower case $\psi_L$ (left) and the upper one $\psi_U$ (right).}\label{fig-unimod}
\includegraphics[width=5cm]{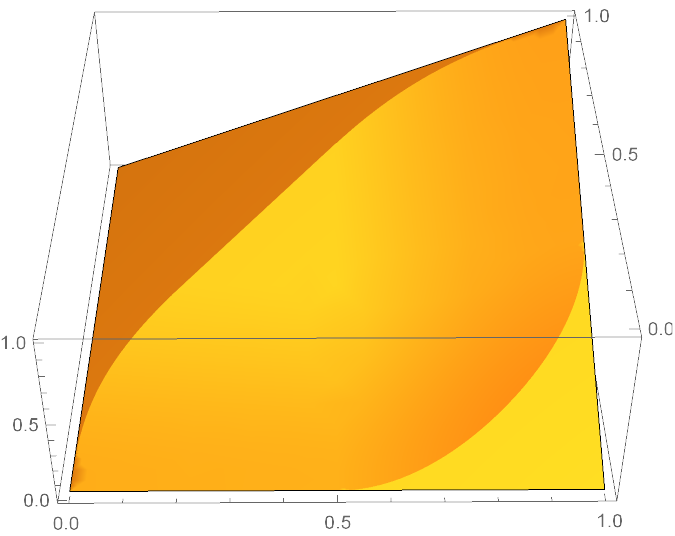}\ \ 
\includegraphics[width=5cm]{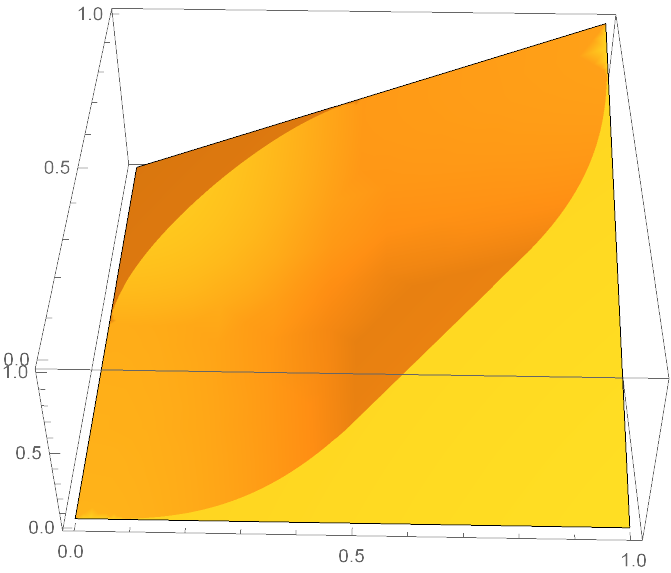}\ \ \includegraphics[width=4cm]{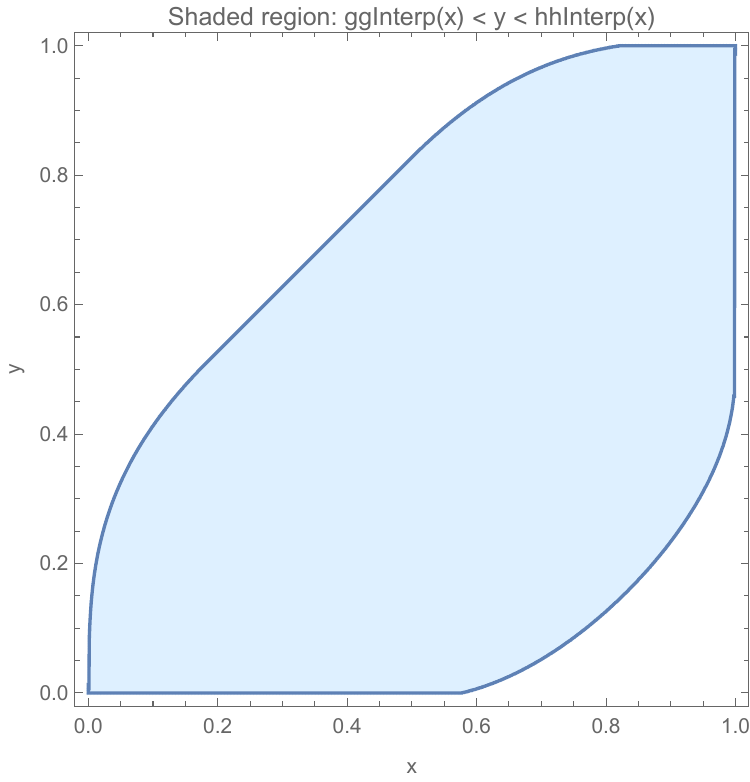}\  \ \ \ \ \ \includegraphics[width=4cm]{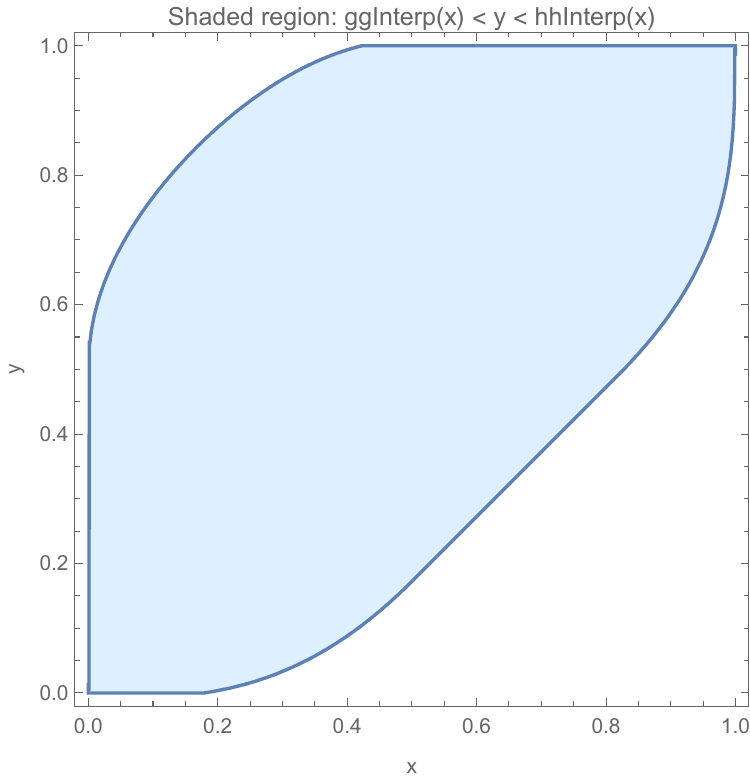}
\label{fig:figure1}
\end{figure}
\end{center}

This result has important consequences on the shape of region $\mathcal{R}_\psi$. Figure \ref{fig-bimod} illustrates the situation when there is exactly one point $x\in(0,1)$ satisfying the conditions of Corollary \ref{cor-region}. This point splits the region into two connected parts. On the other hand, if there is no such point, the region is connected as visualized in Figure \ref{fig-unimod}.

\begin{corollary}\label{cor-region}
  For a point $x\in(0,1)$ the following are equivalent
  \begin{enumerate}[(a)]
    \item $\zeta(x)=0$,
    \item $g(x)=x$,
    \item $h(x)=x$.
  \end{enumerate}
  In this case we have that $C_\psi$ equals function $\min$ on the rectangles $\{(u,v)\mid u\leq x\leq v\}$ and $\{(u,v)\mid v\leq x\leq u\}$. 
\end{corollary}

This corollary is an immediate consequence of Corollary \ref{cor-delta-t-min} and Proposition \ref{prop-g-h}.

It is time to summarize our results. 

\begin{proposition}\label{prop-result}
  \begin{enumerate}[(a)]
    \item The boundary functions are
        \begin{align*}
	   g(x) &= \max\{ y \mid \psi(y) + \zeta(y) \leq \psi(x) \}, \\
	   h(x) &= \max\{ y \mid \xi(y) - \zeta(y) \leq \xi(x) \}, 
        \end{align*}
         and we have \(g(x) \leq x \leq h(x)\) for all $x\in[0,1]$. The inequalities are strict if and only if \(\zeta(x) > 0\).
    \item \begin{align*}
	S_\psi(x, y) &= \begin{cases}
		\psi(x), & \text{if } y > g(x), \\
		\psi(y) + \zeta(y), & \text{if } y \leq g(x),
	\end{cases} \\
	T_\psi(x, y) &= \begin{cases}
		\xi(x), & \text{if } y > h(x), \\
		\xi(y) - \zeta(y), & \text{if } y \leq h(x),
	\end{cases}
\end{align*}
with continuity at \(y = g(x)\) and \(y = h(x)\). 
    \item \begin{align*}\label{eq-C-psi}
	C_\psi(x, y) &= \begin{cases}
		\psi(x) - \psi(y) + \delta(y), & \text{on } \mathcal{R}_\psi 
\\
		M(x,y), & \text{otherwise},
	\end{cases} \\
	&= \min\{ x, y, \psi(x) - \psi(y) + \delta(y) \} \\
	&= \min\{ x, y, \psi(x) - \psi(y) - \zeta(y) + y \}.
\end{align*}
  \end{enumerate} 
\end{proposition}

\begin{proof}
  \emph{(a)} This follows easily from Proposition \ref{prop-g-h}\emph{(a)},\emph{(b)} and Corollary \ref{cor-region}. The functions \(S_\psi\) and \(T_\psi\) from the previous section simplify into
\begin{align*}
	S_\psi(x, y) &= \psi(x) \wedge (\psi(y) + \zeta(y)), \\ 
	T_\psi(x, y) &= \xi(x) \wedge (\xi(y) - \zeta(y)), 
\end{align*}
and \emph{(b)} follows. Corollary~\ref{cor-c-psi-expression-track} provides an explicit expression for the copula \(C_\psi\) implying \emph{(c)}.
\end{proof}



Here is the identity track version of the first of our main results.


\begin{theorem}\label{thm-low-up}
    Let $\delta$ be a function satisfying conditions (a),(b),(c) of a $\varphi$-diagonal with $\varphi(x) = x$ and let $\psi: [0, 1] \to \mathbb{R}$ satisfy $\psi(0) = 0$. Then the following are equivalent: 
	\begin{enumerate}[(a)]
		\item The canonical functions are all increasing.
		\item For all $0\le x\le y\le 1$, we have
		\begin{equation*}
			V^-_{x, y}(\zeta) \leq \psi(y) - \psi(x) \leq (y - x) - V^+_{x, y}(\zeta),
		\end{equation*}
		where $\zeta(x) = x - \delta(x)$.
		\item $C_\psi$ is a copula with diagonal section $\delta$ corresponding to track $\varphi(x) = x$. 
	\end{enumerate}	
\end{theorem}

In the next section we give the general track version of this theorem, so that it will become a simple corollary of Theorem \ref{thm-bounds-track} and we may omit its proof now.
We see that there are certain bounds on $\psi$ depending on $\delta$ (or equivalently $\zeta$) in order to consider it for the desired role. Denote by \(\low\psi\) and \(\up\psi\) the respective lower and upper bound for $\psi$. It is clear that not every possible function between these two bounds \mo{is eligible}. 

\begin{corollary}
	We have: 
	\begin{align*}
		\low\psi(x) &= V^+_{0, x}(-\zeta) = V^-_{0, x}(\zeta), \\
		\up\psi(x) &= x - V^+_{0, x}(\zeta).
	\end{align*}
    For any copula \(C\) with \(C(x, x) = \delta(x)\):
	\begin{enumerate}[(a)]
		\item \(C_{\up\psi}(x, y) \geq C(x, y)\) if \(y \geq x\),
		\item \(C_{\low\psi}(x, y) \geq C(x, y)\) if \(y \leq x\).
	\end{enumerate}
\end{corollary}

The proof of this corollary will also be omitted since it is an immediate consequence of its general track version Corollary \ref{cor-track-region} given in the next section.

Define \(\kappa(x, y) = \psi(x) - \psi(y) - \zeta(y) + y\). Since \(\kappa\) increases with \(\psi(x) - \psi(y)\), it is maximized by \(\up\psi\) for \(x \geq y\) and \(\low\psi\) for \(x \leq y\). In the next two displayed computations we cite equations \emph{(f)} and \emph{(g)} stated just before Lemma \ref{lem-posvar},

For \(x \leq y\) we have:
\begin{align*}
	\kappa(x, y) &= \low\psi(x) - \low\psi(y) - \zeta(y) + y \\
	&= V^-_{0, x}(\zeta) - V^-_{0, y}(\zeta) - \zeta(y) + y \\
	&= -V^-_{x, y}(\zeta) + y - \zeta(y) \\
	&= y - \zeta(y) - \frac{1}{2} (TV_{x, y}(\zeta) + \zeta(x) - \zeta(y)) \quad (\text{by \emph{(g)}}) \\
	&= y - \frac{1}{2} (TV_{x, y}(\zeta) + \zeta(x) + \zeta(y)).
\end{align*}

For \(x \geq y\) we have:
\begin{align*}
	\kappa(x, y) &= \up\psi(x) - \up\psi(y) - \zeta(y) + y \\
	&= x - V^+_{0, x}(\zeta) - (y - V^+_{0, y}(\zeta)) - \zeta(y) + y \\
	&= x - V^+_{y, x}(\zeta) - \zeta(y) \\
	&= x - \frac{1}{2} (TV_{y, x}(\zeta) + \zeta(x) - \zeta(y)) - \zeta(y) \quad (\text{by \emph{(f)}}) \\
	&= x - \frac{1}{2} (TV_{y, x}(\zeta) + \zeta(x) + \zeta(y)).
\end{align*}

Thus, the upper bound is attained by \(C_{\low\psi}(x, y)\) for \(x \leq y\) and \(C_{\up\psi}(x, y)\) for \(x \geq y\), as \(C_\psi(x, y) = \min\{ x, y, \kappa(x, y) \}\) is maximized accordingly. One should compare this with the result obtained in \cite{KoBuMoSt}, where the lower and upper bounds of (quasi-)copulas with given diagonal section are studied. While they only get the extremal undominated quasi-copulas in \cite[Equation (3)]{KoBuMoSt}, which is not surprisingly the same as ours just obtained, we are aimed at getting all (!) undominated copulas and consequently quasi-copulas. 

As a side result we get a reformulation of a known fact.

\begin{corollary}
	Let \(\delta: [0, 1] \to [0, 1]\) be an increasing function, and set \(\zeta(x) = x - \delta(x)\). A copula \(C\) with \(\delta(x) = C(x, x)\)  exists if and only if:
	\[
	TV_{x, y}(\zeta) \leq y - x \quad \text{for all } 0 \leq x \leq y \leq 1.
	\]
	This is equivalent to the property that  \(\zeta\) is 1-Lipschitz, which is also equivalent to the property that  \(\delta\) is 2-Lipschitz,
\end{corollary}
\begin{proof}
	We first show the equivalence of total variation Lipschitz conditions. It is clear that the condition \(TV_{x, y}(\zeta) \leq y - x\) for all \(x \leq y\) is equivalent to \(\zeta\) being 1-Lipschitz. (Moreover, an easy application of the triangle inequality shows that this is equivalent to $\delta$ being 2-Lipschitz.) This follows from the fact that the total variation of \(\zeta\) over an interval is bounded by the difference between the endpoints of the interval, which is the defining property of a 1-Lipschitz function.
	
	Further, a copula \(C_\psi\) with diagonal \(\delta\) exists if and only if there exists \(\psi\) such that the following functions
	\[
	\psi(x), \quad x - \psi(x), \quad \psi(x) + \zeta(x), \quad x - \psi(x) - \zeta(x),
	\]
	 are increasing; \mo{or equivalently, such that the bounds} 
	\[
	V^-_{x, y}(\zeta) \leq \psi(y) - \psi(x) \leq (y - x) - V^+_{x, y}(\zeta),
	\]
	hold for all \(x \leq y\), as stated in Theorem~\ref{thm-low-up}. 
	
	For such \(\psi\) to exist, the interval must be non-empty, i.e., we must have
	\[
	V^-_{x, y}(\zeta) \leq (y - x) - V^+_{x, y}(\zeta),
	\]
	which is equivalent to
	\[
	V^-_{x, y}(\zeta) + V^+_{x, y}(\zeta) = TV_{x, y}(\zeta) \leq y - x,
	\]
	using the identity \(V^-_{x, y}(\zeta) + V^+_{x, y}(\zeta) = TV_{x, y}(\zeta)\) from variational property \emph{(a)} of Section \ref{s-prlm}.
\end{proof}

\section{All undominated copulas with a given diagonal section}\label{sec-final}

In this section we present all possible undominated copulas with a given diagonal section. We have seen  that every undominated copula having the same diagonal section is of the form $C_\psi$. We will now show that every copula of the form $C_\psi$ for an appropriate function $\psi$ is undominated. This will only be shown for the standard track $\varphi(x)=x$. However, we first develop some additional tools for general tracks.

Fix a track section \(\delta\) on \(B_\varphi = \{(x, \varphi(x)) \mid x \in [0, 1]\}\), where \(\varphi: [0, 1] \to [0, 1]\) is strictly increasing, and let $C \in \mathcal{C}(B_\varphi, \delta)$. Note that \(y - \delta(\varphi^{-1}(y)) + \psi(\varphi^{-1}(y))\) is increasing in \(y\) if and only if \(\varphi(x) - \delta(x) + \psi(x)\) is increasing, since \(y = \varphi(x)\) and \(\varphi^{-1}\) is increasing. Similarly, \(\delta(\varphi^{-1}(y)) - \psi(\varphi^{-1}(y))\) is increasing if and only if \(\delta(x) - \psi(x)\) is increasing.

\begin{theorem}\label{thm-bounds-track}
	Let $\delta$ be a function satisfying conditions (a),(b),(c) of a $\varphi$-diagonal and let $\psi: [0, 1] \to \mathbb{R}$ satisfy $\psi(0) = 0$. Then the following are equivalent: 
	\begin{enumerate}[(a)]
		\item The canonical functions are all increasing.
		\item For all $0\le x\le y\le 1$,
		\begin{equation*}
			V^-_{x, y}(\varphi - \delta) \leq \psi(y) - \psi(x) \leq (y - x) - V^+_{x, y}(\zeta),
		\end{equation*}
		where $\zeta(x) = x - \delta(x)$.
		\item $C_\psi$ is a copula with diagonal section $\delta$ corresponding to track $\varphi$. 
	\end{enumerate}	
\end{theorem}

\begin{proof}
	We start with the implication \emph{(a)} $\implies$ \emph{(b)}. Since $\psi$ and $\chi = \psi - (\delta - \varphi) = \psi - \delta + \varphi$ are increasing, Lemma~\ref{lem-posvar} gives:
	\[
	\psi(y) - \psi(x) \geq V^+_{x, y}(\delta-\varphi) = V^-_{x, y} (\varphi-\delta).
	\]
	
	Similarly, since $\xi = x - \psi$ and $\eta = \delta -\psi = \xi - \zeta$ are increasing:
	\[
	\xi(y) - \xi(x) = (y - \psi(y)) - (x - \psi(x)) \geq V^+_{x, y}(\zeta),
	\]
	so $\psi(y) - \psi(x) \leq (y - x) - V^+_{x, y}(\zeta)$.
	
	Thus:
	\[
	V^-_{x, y}(\varphi-\delta) \leq \psi(y) - \psi(x) \leq (y - x) - V^+_{x, y}(\zeta),
	\]
	confirming \emph{(b)}.
	
	We prove the reverse implication \emph{(b)} $\implies$ \emph{(a)} by proving monotonicity of each canonical function separately. 
\begin{itemize}
  \item[$-$] {$\psi$ is increasing:} The left inequality in \emph{(b)} 
  gives $\psi(y)-\psi(x) \ge V^-_{x, y}(\varphi - \delta) \ge 0$ for all $0\le x\le y \le 1$.   
  \item[$-$] {$\eta = \delta - \psi$ is increasing:} Take some $0\le x\le y \le 1$. First note that $\delta(y)-\delta(x) = (y-x)-(\zeta(y)-\zeta(x))$. It is also clear that $V^+_{x, y}(\zeta) \ge \zeta(y)-\zeta(x)$. Hence, $(y-x)-V^+(\zeta) \le (y-x) - (\zeta(y)-\zeta(x)) = \delta(y)-\delta(x)$. By the second inequality in \emph{(b)} 
      we have that $\psi(y)-\psi(x) \le \delta(y)-\delta(x)$. Therefore, $\eta(y)-\eta(x) = (\delta(y)-\psi(y)) - (\delta(x)-\psi(x)) \ge 0$. 
  \item[$-$] {$\xi(x) = x - \psi(x)$ is increasing:} By rearranging the second inequality in \emph{(b)} 
      we obtain $\xi(y)-\xi(x) = (y-x)-(\psi(y)-\psi(x)) \ge V^+_{x, y}(\zeta) \ge 0$. 
  \item[$-$] {$\chi(x) = \psi - \delta + \varphi$ is increasing:} This follows directly by applying Lemma~\ref{lem-posvar} on inequality $\psi(y)-\psi(x) \ge V^-_{x, y}(\varphi - \delta)$, which gives exactly the condition that $\chi = \psi - (\delta - \varphi)$ is increasing.
\end{itemize}
	
	The implication \emph{(a)} $\implies$ \emph{(c)} is shown by Theorem~\ref{prop-copula-psi}. The reverse implication  \emph{(c)} $\implies$ \emph{(a)} follows directly from Proposition~\ref{prop-cont-initial}. 
\end{proof}

\begin{corollary}\label{cor-jwaid}

	Let \(\varphi: [0, 1] \to [0, 1]\) be increasing, \(B_\varphi\) the corresponding track, and \(\delta(x)\) a $\varphi$-diagonal. Recall \(\zeta(x) = x - \delta(x)\)
. A copula \(C\) with \(C(x, \varphi(x)) = \delta(x)\) exists if and only if:
	\[
	V^-_{x, y}(\varphi(x) - \delta(x)) + V^+_{x, y}(\zeta) \leq y - x \quad \text{for all } 0 \leq x \leq y \leq 1.
	\]
\end{corollary}
\begin{proof}
	The interval given in Theorem \ref{thm-bounds-track} is non-empty if:
	\[
	V^-_{x, y}(\varphi(x) - \delta(x)) \leq (y - x) - V^+_{x, y}(\zeta),
	\]
	i.e., \(V^-_{x, y}(\tilde{\delta}) + V^+_{x, y}(\zeta) \leq y - x\). This ensures that some \(\psi\) exists, yielding \(C_\psi\) with \(C_\psi(x, \varphi(x)) = \delta(x)\).
\end{proof}

This corollary helps us extend a known result on diagonal sections corresponding to general tracks \cite{DeBaDeMeJw}.

\begin{proposition}\label{prop-track-existence}
	Let \(\varphi: [0, 1] \to [0, 1]\) be strictly increasing and continuous, and \(\delta: [0, 1] \to [0, 1]\) be continuous, with \(\zeta(x) = x - \delta(x)\)
. The following are equivalent:
	\begin{enumerate}[(a)]
		\item A copula \(C\) with \(C(x, \varphi(x)) = \delta(x)\) exists,
		\item \(V^-_{x, y}(\varphi(x) - \delta(x)) + V^+_{x, y}(\zeta) \leq y - x\) for all \(0 \leq x \leq y \leq 1\),
		\item \(\delta(y) - \delta(x) \leq y - x + \varphi(y) - \varphi(x)\) for all \(0 \leq x \leq y \leq 1\).
	\end{enumerate}
\end{proposition}

\begin{proof}
    Corollary \ref{cor-jwaid} gives equivalence of {\emph{(a)}} and \emph{(b)}, \cite[Proposition 2.1]{DeBaDeMeJw} equivalence of \emph{(a)} and \emph{(c)}.
\end{proof}

\begin{corollary}\label{cor-track-region}
	Let \(\varphi: [0, 1] \to [0, 1]\) be strictly increasing, and \(\delta\) satisfy Proposition~\ref{prop-track-existence}. Define:
	\[
	\low\psi(x) = V^-_{0, x}(\varphi(x) - \delta(x)), \quad \up\psi(x) = x - V^+_{0, x}(\zeta),
	\]
	where \(\zeta(x) = x - \delta(x)\), \(\tilde{\delta}(x) = \varphi(x) - \delta(x)\). For any copula \(C\) with \(C(x, \varphi(x)) = \delta(x)\):
	\begin{enumerate}[(a)]
		\item \(C_{\up\psi}(x, y) \geq C(x, y)\) if \(y \geq \varphi(x)\),
		\item \(C_{\low\psi}(x, y) \geq C(x, y)\) if \(y \leq \varphi(x)\).
	\end{enumerate}
\end{corollary}
\begin{proof}
	By Theorem~\ref{thm-bounds-track}, for \(x \leq y\):
	\[
	V^+_{x, y}(\varphi - \delta) \leq \psi(y) - \psi(x) \leq (y - x) - V^+_{x, y}(\zeta),
	\]
\ds{for any eligible function $\psi$, i.e., such that the four corresponding canonical functions are increasing. S}o \(\low\psi(y) - \low\psi(x) = V^-_{x, y}(\tilde{\delta})\) is minimal, and \(\up\psi(y) - \up\psi(x) = (y - x) - V^+_{x, y}(\zeta)\) is maximal among such \(\psi\).
	
	Since \(C_\psi(x, y) = \min\{x, y, \psi(x) - \psi(\varphi^{-1}(y)) + \delta(\varphi^{-1}(y))\}\) we have the following. 
\mo{For \(y \geq \varphi(x)\), we have \(\varphi^{-1}(y) \geq x,\) \[\up\psi(x) - \up\psi(\varphi^{-1}(y)) \geq \psi(x) - \psi(\varphi^{-1}(y)),\] so that \(C_{\up\psi}(x, y) \geq C_\psi(x, y)\).
	For \(y \leq \varphi(x)\) we have \(\varphi^{-1}(y) \leq x,\)\[ \low\psi(x) - \low\psi(\varphi^{-1}(y)) \geq \psi(x) - \psi(\varphi^{-1}(y)),\] so that \(C_{\low\psi}(x, y) \geq C_\psi(x, y)\).}
	
	If \(y = \varphi(x)\), \(C_\psi(x, y) = \delta(x)\), consistent for all \(\psi\).
	
	Since for any copula $C$, by Corollary~\ref{cor-c-psi-undominated}, there exists some $C_\psi \ge C$, it follows that the upper bounds for copulas of the form $C_\psi$ are also upper bounds for all copulas. 
\end{proof}

\ds{Let us denote the set of all eligible functions $\psi$ by $\mathcal{E}$. Then it is ordered naturally pointwise with respect to the only variable. Functions $\low\psi$ and $\up\psi$ are respectively the minimal and the maximal point of $\mathcal{E}$ in this order. Observe that this order induces also an order on the set of copulas $\{C_\psi\}_{\psi\in\mathcal{E}}$. }
\smallskip


It is time to give our main result, this time for the case of identity track only. But let us start with a technical lemma.

\begin{lemma}\label{lem-maximality}
	Let \(\psi_1, \psi_2: [0, 1] \to \mathbb{R}\) satisfy the conclusion of Proposition~\ref{prop-cont-initial} (i.e., \(\psi_i(x)\), \(x - \psi_i(x)\), \(\psi_i(x) + \zeta(x)\), \(\delta(x) - \psi_i(x)\) are increasing for $i=1,2$), and let \([a, b] \subseteq [0, 1]\) be an interval where \(\delta(x) < x - c\) for some \(c > 0\). If \(\psi_1(t) - \psi_1(a) \neq \psi_2(t) - \psi_2(a)\) for some \(t \in (a, b)\), then:
	\begin{enumerate}[(a)]
		\item \(C_{\psi_1} \neq C_{\psi_2}\),
		\item neither \(C_{\psi_1} > C_{\psi_2}\) nor \(C_{\psi_2} > C_{\psi_1}\) holds everywhere,
		\item there exist \(u, v \in [0, 1]\) such that:
		\begin{equation*}
			(C_{\psi_1}(u, v) - C_{\psi_2}(u, v))(C_{\psi_1}(v, u) - C_{\psi_2}(v, u)) < 0.
		\end{equation*}
	\end{enumerate}
\end{lemma}
\begin{proof}
	Define \(g_i(x) = \max\{y \mid \psi_i(y) + \zeta(y) \leq \psi_i(x)\}\), \(h_i(x) = \min\{y \mid \delta(y) - \psi_i(y) \leq x - \psi_i(x)\}\). Since \(\zeta(x) = x - \delta(x) > c\),
	\(\psi_i(x) + \zeta(x) > \psi_i(x) + c > \psi_i(x)\), so \(g_i(x) < x\), and 
	\(\delta(x) - \psi_i(x) < x - c - \psi_i(x) < x - \psi_i(x)\), so \(h_i(x) > x\).
	Let \(\alpha = \min_{i=1,2} \min_{x \in [a, b]} \{(x - g_i(x)) \wedge (h_i(x) - x)\} > 0\).
	
	Construct a partition \(a = x_0 < x_1 < \cdots < x_n = b\) with \(x_k - x_{k-1} < \alpha\), and \(t = x_k\) for some \(k\). Since \(\psi_1(t) - \psi_1(a) \neq \psi_2(t) - \psi_2(a)\), there exists \(k\) such that \(\psi_1(x_k) - \psi_1(x_{k-1}) \neq \psi_2(x_k) - \psi_2(x_{k-1})\). For \(u = x_{k-1}\), \(v = x_k\),  \(v < h_i(u)\), so \(C_{\psi_i}(u, v) = \psi_i(u) - \psi_i(v) + \delta(v)\), and \(u < g_i(v)\), so \(C_{\psi_i}(v, u) = \psi_i(v) - \psi_i(u) + \delta(u)\).
	If \(\psi_1(v) - \psi_1(u) < \psi_2(v) - \psi_2(u)\), then \(C_{\psi_1}(v, u) < C_{\psi_2}(v, u)\), \(C_{\psi_1}(u, v) > C_{\psi_2}(u, v)\), and the product is negative.
	Thus, \(C_{\psi_1} \neq C_{\psi_2}\), and no uniform ordering exists.
\end{proof}

\begin{theorem}\label{thm-undominated}
	Let \(C_{\psi_1}\) and \(C_{\psi_2}\) be copulas as in Lemma~\ref{lem-maximality}. Then either \(C_{\psi_1} = C_{\psi_2}\) or 
neither dominates the other.
\end{theorem}
\begin{proof}
	If \(C_{\psi_1}(x, y) \neq C_{\psi_2}(x, y)\), assume \(x < y\) (since \(C_{\psi_i}(x, x) = \delta(x)\)) and \(C_{\psi_1}(x, y) < C_{\psi_2}(x, y) \leq x\). By continuity of \(\zeta\) and Corollary~\ref{cor-delta-t-min}, there exists \([a, b] \supseteq [x, y]\) with \(\zeta(t) > c > 0\). Since \(C_{\psi_1}(x, y) < x\),  \(C_{\psi_1}(x, y) = \psi_1(x) - \psi_1(y) + \delta(y)\), and \(C_{\psi_2}(x, y) = \psi_2(x) - \psi_2(y) + \delta(y) > C_{\psi_1}(x, y)\). So \(\psi_1(x) - \psi_1(y) < \psi_2(x) - \psi_2(y)\).
	Applying Lemma~\ref{lem-maximality} with \(a = x\), \(t = y\), there exist \(u, v\) satisfying 
condition \emph{(c)} of the Lemma.
\end{proof}

\ds{
\begin{corollary}\label{cor-equiv}
    In case of identity track a copula with a given diagonal section $\delta$ is undominated in the family of copulas with the same $\delta$ if and only if it is of the form $C_\psi$ for some $\psi\in\mathcal{E}$.
\end{corollary}

\begin{proof}
  By Corollary \ref{cor-c-psi-undominated} every undominated copula is of the form $C_\psi$ for some $\psi\in\mathcal{E}$. Conversely, choose any copula $C_{\psi_1}$ with a given diagonal section $\delta$ and $\psi_1\in\mathcal{E}$. Adjoin the canonical quadruple to copula $C=C_{\psi_1}$ and denote the first one of the four canonical functions by $\psi_2$. So, $C_{\psi_2}\ge C=C_{\psi_1}$ is undominated. However, this can only be true when $C_{\psi_2}=C_{\psi_1}$ by Theorem \ref{thm-undominated} and consequently $C_{\psi_1}$ is undominated.
\end{proof}

Finally, we give a brief insight into the possible study of undominated quasi-copulas. Assume $\mathcal{E}$ is the eligible set of functions corresponding to a track function $\varphi$ and $\varphi$-diagonal $\delta$. Moreover, choose two eligible functions $\psi_U,\psi_L\in\mathcal{E}$. Define diagonal splice of the corresponding copulas \cite{NeQuMoRoLaUbFl} (cf. also \cite{JwDeMeHaIsDeBa})
\[
    C_{\psi_L}^{\psi_U}(u,v)=\begin{cases}
                               C_{\psi_U}(u,v), & \mbox{if $v\ge\varphi(u)$}  \\
                               C_{\psi_L}(u,v), & \mbox{otherwise}.
                             \end{cases}
\]
It is not hard to see (and is also shown in \cite{NeQuMoRoLaUbFl}) that this is always a quasi-copula.
A natural question is, whether it is undominated.
}


\begin{thebibliography}{00}

\bibitem{Bert} S.\ Bertino, On dissimilarity between cyclic permutations, Metron \textbf{35} (1–2) (1977) 53--88, MR600402.

\bibitem{DeBaDeMeJw} B.\ De Baets, H.\ De Meyer, T.\ Jwaid, On the degree of asymmetry of a quasi-copula with respect to a curve, Fuzzy Sets and Systems \textbf{354} (2019) 84–103.


\bibitem{DeDoUbFl} A.\ Dehgani, A.\ Dolati, M.\ \'{U}beda-Flores, \emph{Measures of radial asymmetry for bivariate random vectors}, Stat.\ Papers \textbf{54} (2013), 271--286.
    
\bibitem{DuKoMeSe} F.~Durante, A.~Koles\' arov\' a, R.~Mesiar, C.~Sempi, {\emph{Copulas with given diagonal sections: novel constructions and applications}}, International Journal of Uncertainty, Fuzziness and Knowledge-Based Systems, (2007), 14 pp.


\bibitem{DuSe} F.~Durante, C.~Sempi, Principles of Copula Theory, CRC/Chapman \& Hall, Boca Raton (2015).

\bibitem{FeSaUbFl} J. Fern\' andez-S\' anchez, M. \' Ubeda-Flores, \emph{Copulas with given track and opposite track sections: Solution to a problem on diagonals}, Fuzzy Sets and Systems \textbf{308} (2017), 133–137

\bibitem{feSaUbFl1} J.\ Fern\' andez-S\' anchez, M.\ \' Ubeda Flores, Constructions of copulas with given diagonal (and opposite diagonal) sections and some generalizations, Depend.\ Model. \textbf{6 }(1) (2018) 139--155, MR3824791.
    
\bibitem{FrNe} G.\ A.\ Fredricks, R.\ B.\ Nelsen, Copulas constructed from diagonal sections, in: Distributions with Given Marginals and Moment Problems, Prague, 1996, Kluwer Acad. Publ., Dordrecht, 1997, 129--136, MR1614666.

\bibitem{FrNe1} G.\ A.\ Fredricks, R.\ B.\ Nelsen, The Bertino family of copulas, in: Distributions with Given Marginals and Statistical Modelling, Kluwer Acad. Publ., Dordrecht, 2002, pp. 81–91, MR2058982.

\bibitem{FuSc} S.~Fuchs, K.~D.~Schmidt, \textsl{Bivariate copulas, transformations, asymmetry and measures of concordance}, Kybernetika \textbf{50} (2014), 109--125

\bibitem{GeNe} C.\ Genest, J.\ Ne\v{s}lehov\'a, \textsl{Assessing and Modeling Asymmetry in Bivariate Continuous Data}, In: P.\ Jaworski, F.\ Durante, W.K.\ H\"ardle, (eds.), {C}opulae in {M}athematical and {Q}uantitative {F}inance, Lecture Notes in Statistics, Springer Berlin Heidelberg, (2013), 152--16891--114.

\bibitem{JwDeMeHaIsDeBa} T. Jwaid, H. De Meyer, A. Haj Ismail, B. De Baets, \emph{Curved splicing of copulas}, Information Sciences 556 (2021) 95–110.

\bibitem{KlMe} E.\ P.\ Klement and R.\ Mesiar. \textsl{How non-symmetric can a copula be?} Comment.\ Math.\ Univ.\ Carolin.\ \textbf{47} (2006), no.\ 1, 141--148.
    
    
\bibitem{KoBuMoSt} D.\ Kokol Bukovšek, B.\ Mojškerc, N. Stopar, \textsl{Exact upper bound for copulas with a given diagonal section}, Fuzzy Sets and Systems \textbf{480} (2024) 108865.



\bibitem{Loja} S.~\L ojasiewicz, An Introduction to the Theory of Real Functions, John Wiley \& Sons (1988).








\bibitem{Nels} R.\ B.\ Nelsen, An introduction to copulas, 2nd edition, Springer-Verlag, New York (2006).

\bibitem{Nels1} R.\ B.\ Nelsen. {\it Extremes of nonexchangeability.} Statist. Papers {\bf 48} (2007), no. 2, 329--336.

\bibitem{NeQuMoRoLaUbFl} R.\ B.\ Nelsen, J.\ J.\ Quesada-Molina, J.\ A.\ Rodr\' iguez-Lallena, M.\ \' Ubeda-Flores, On the construction of copulas and quasi-copulas with given diagonal sections, Insurance: Mathematics and Economics \textbf{42} (2008) 473--483.




\bibitem{OmSt4} M.\ Omladi\v{c}, N.\ Stopar, \emph{Multivariate imprecise Sklar type theorems}, Fuzzy Sets and Systems, \textbf{428} (2022) 80--101.


\bibitem{OuSuZh} Y.\ Ouyang, Y.\ Sun, H.-P.\ Zhang, On an upper bound of the set of copulas with a given curvilinear section, Fuzzy Sets and Systems 500 (2025) 109199



\bibitem{Skla} A.\ Sklar, Fonctions de r\'{e}partition \`{a} $n$ dimensions et leurs marges, Publ.\ Inst.\ Stat.\ Univ.\ Paris \textbf{8} (1959) 229--231.

\bibitem{ZoSuXi} W. Zou, L. Sun, J. Xie, Best-possible bounds on the sets of copulas and quasi-copulas with given curvilinear sections, Fuzzy Sets Syst. 441 (2022) 335–365.

\bibitem{math} Wolfram Research, Inc. Mathematica, Version {\bf 11}, Champaign, IL, 2017.

\end{thebibliography}
\end{document}